%
%
%
\documentstyle{amsppt}

\magnification=1200
\NoBlackBoxes
\def\<{\langle}
\def\>{\rangle}
\def\L{\bold L}
\def\M{\bold M}
\def\P{{\Bbb P}}
\def\N{{\Bbb N}}

\topmatter
\title Orthogonal matrix polynomials and higher order
recurrence relations \endtitle
\author A.J. Dur\'an and W. Van Assche
\endauthor
\affil Universidad de Sevilla \\ and \\
  Katholieke Universiteit Leuven
\endaffil
\thanks
W.V.A. is a Senior Research Associate of the Belgian National Fund for
Scientific Research
\endthanks
\address (A.J.D.) Departamento de An\'alisis Matem\'atico, Universidad de
Sevilla, Apdo.\ 1160, E--41080 Sevilla, SPAIN
\endaddress
\email duran\@ cica.es \endemail
\address (W.V.A.) Department of Mathematics, Katholieke Universiteit
Leuven, Celestijnenlaan 200 B, B--3001 Heverlee (Leuven), BELGIUM
\endaddress
\email fgaee03\@ cc1.KULeuven.ac.be \endemail
\rightheadtext{Matrix polynomials and recurrence relations}
\abstract
It is well-known that orthogonal polynomials on the real line
satisfy a three-term recurrence relation and conversely
every system of polynomials satisfying a three-term
recurrence relation is orthogonal with respect to some
positive Borel measure on the real line. In this paper we
extend this result and show that every system of
polynomials satisfying some $(2N+1)$-term
recurrence relation can be expressed in terms of orthonormal
matrix polynomials for which the coefficients are $N\times
N$ matrices. We apply this result to polynomials
orthogonal with respect to a discrete Sobolev inner product
and other inner products in the linear space of polynomials.
As an application we give a short proof of Krein's characterization
of orthogonal polynomials with a spectrum having a finite number
of accumulation points.
\endabstract

\endtopmatter

\document

\head 1. Introduction \endhead

A sequence of orthonormal polynomials $p_n(x)$
$(n=0,1,2,\ldots)$ on the real line, orthonormal with some probability
measure $\mu$, always satisfies a three-term recurrence relation
$$   xp_n(x) = a_{n+1} p_{n+1}(x) + b_n p_n(x) + a_{n} p_{n-1}(x), \tag 1.1 $$
with initial conditions $p_{-1}(x) = 0$ and $p_0(x) = 1$. The recurrence
coefficients are given by
$$   a_n = \int x p_{n-1}(x)p_n(x)\, d\mu(x) > 0, \quad
   b_n = \int x p_n^2(x) \, d\mu(x) \in {\Bbb R} . $$
The converse is also true: a system of polynomials satisfying a three-term
recurrence relation \thetag{1.1} with $a_{n+1} > 0$ and $b_n \in {\Bbb R}$
$(n=0,1,2,\ldots)$ is always a system of orthonormal polynomials with
respect to some probability measure $\mu$ on the real line. This converse
result was given by Favard in 1935 \cite{10} but was known earlier and
appears already in the books by Stone \cite{21, Theorem 10.27 on pp.\
545--546}, Perron \cite{18, \S 36, Satz 4.6/4.7} and Wintner \cite{22, \S 32
and \S 87} who attributes the case of finite support to E. Heine \cite{12, \S
108}. In an elementary form (lacking the Riesz
representation theorem) the result is in Stieltjes' 1894 work \cite{20, \S 11}
and apparently the result was already known by Chebyshev \cite{5} when the
support of the measure $\mu$ is finite.

The importance of Favard's theorem is that orthogonal polynomials and
polynomials satisfying a three-term recurrence relation are the same thing.
Properties regarding zeros (real and simple zeros, interlacing of zeros) and
positivity of connection coefficients can thus be studied from two points of
view: on one hand using the orthogonality, on the other hand using the
recurrence relation.

Favard's theorem can also be proved for orthogonal matrix polynomials.
Orthogonal matrix polynomials on the real line have been considered
in detail by M.G. Krein \cite{13}. See also the book by Berezanski\u\i{}
\cite{4} and more recent papers by
Aptekarev and Nikishin \cite{1}, Geronimo \cite{11}, and
Sinap and Van Assche \cite{19}.
Consider  matrix polynomials
satisfying the three-term recurrence relation
$$   x P_k(x) = D_{k+1} P_{k+1}(x) + E_k P_k(x) + D_k^* P_{k-1}(x), \tag 1.2
$$ with $P_0(x)=I$ and $P_{-1}(x)= 0$,
where $P_k(x)$ are matrix polynomials with coefficients in ${\Bbb C}^{N\times
N}$ and the recurrence coefficients $D_{k+1}, E_k$ are also
$N\times N$ matrices for which $E_k^*=E_k$ and $\det D_k \neq 0$.
By using  spectral theory Aptekarev and Nikishin \cite{1} show that the
polynomials $P_n(x)$ $(n=0,1,2,\ldots)$ are orthonormal with respect to some
Hermitian matrix of measures $\M = (\mu_{k,l})_{k,l=1}^N$ which is positive
definite:
$$   \int P_n(x)\, d\M(x) \, P_m^*(x) = \delta_{m,n} I . $$

Recently one of us studied polynomials $p_n(x)$ satisfying a $(2N+1)$-term
recurrence relation
$$ h(x)  p_n(x) = c_{n,0} p_n(x) + \sum_{k=1}^N \left(
  c_{n,k} p_{n-k}(x) + c_{n+k,k} p_{n+k}(x) \right) , \tag 1.3  $$
  where $h$ is a polynomial of degree $N$
and $c_{n,k}$ $(n=0,1,2,\ldots)$
are real sequences for $k=0,1,\ldots,N$ with $c_{n,N} \neq 0$. From \cite{7}
we get, after straightforward reformulation, that  a sequence of polynomials
satisfies a $(2N+1)$-term recurrence relation if and only if
the following orthogonality
condition holds:
there exists a $N\times N$ matrix of measures
$\mu =( \mu_{k,l})_{k,l=0}^{N-1}$ such that
the bilinear form
$$
B_\mu (p,q)=\sum _{k,l=0}^{N-1}\int R_{h,N,k}(p)\overline{R_{h,N,l}(q)}\,d\mu
_{k,l}, \tag 1.4
$$
where
$$  R_{h,N,k} (p)(x) = \sum_{n=0}^m a_{k,n} x^n \quad  \text {if} \quad
    p(x) = \sum_{k=0}^{N-1} \sum_{n=0}^m a_{k,n} x^k h^n(x) \tag 1.5 $$
 is an inner product on the linear space of polynomials
 $\P$ and $(p_n)_n$ is the
sequence of orthonormal polynomials with respect to $B_\mu $.
 In \cite{8}          this is improved showing that the matrix of
 measures $\mu $ can be taken to be positive definite.

In Section 2 we will show that polynomials satisfying a higher
order recurrence relation of the form \thetag{1.3} (but with complex
coefficients) are closely
related to matrix polynomials satisfying a three-term recurrence relation
 and that Favard's theorem for matrix polynomials and Favard's
theorem for polynomials satisfying a higher order recurrence relation are
the same.

In Section
3 we show how discrete Sobolev orthogonal polynomials and some new
orthogonal polynomials with respect to an inner product of the form
$$  \langle f,g \rangle = \int f(x)g(x)\, d\mu(x)
   + \left( \sum_{k=1}^N a_k f(c_k) \right) \left( \sum_{k=1}^N a_k g(c_k)
\right), \qquad a_k, c_k \in {\Bbb R} $$
are related to orthogonal
matrix polynomials and we explicitly give the orthogonality matrix of measures
$\M$ in terms of the parameters in the inner product. The advantage of working
with
the matrix  polynomials is that the orthogonality conditions no longer require
the evaluation of a function and its derivatives in various points. Finally, in
Section
4 we show how Krein's theorem regarding orthogonal polynomials with a spectrum
with finitely many accumulation points follows easily by this correspondence
between higher order recurrence relations and orthogonal matrix polynomials
on the real line.

The connection between matrix polynomials and scalar polynomials obtained
in the present paper shows some analogy  with the
connection between orthogonal polynomials on a lemniscate and orthogonal
polynomials on the unit circle, as given by Marcell\'an and Rodr\'\i guez
\cite{15} and the connection between orthogonal polynomials on an
algebraic harmonic curve and orthogonal matrix polynomials on the real line,
as given by Marcell\'an and Sansigre \cite{17}. In fact the particular choice
of basis in the linear space of polynomials that we use in Sections 2 and 3
is the same as the basis used in \cite{15} and \cite{17}, but the
polynomial $h(x)$ in the present paper is connected with the recurrence
relation, whereas in \cite{15} this polynomial describes the lemniscate
and in \cite{17} it describes the algebraic harmonic curve.

\head 2. Recurrence relation and matrix polynomials
\endhead

In order to establish the main theorem we need to consider the operators
$R_{h,N,k}$
defined in the introduction of this paper (see \thetag{1.5}). For the sake of
simplicity
we start by considering the case $h(x)=x^N$. Then, we denote the operators by
$R_{N,m}$  and it is not hard to see that they are defined by
$$
R_{N,m}(p)(x)=\sum _n{p^{(nN+m)}(0)\over (nN+m)!}x^n ,
$$
i.e., the operator $R_{N,m}$ takes from $p$ just those powers with remainder
$m$ modulo $N$ and then, removes $x^m$ and changes $x^N$ to $x$.
Thus, we have
$$
p(x)=R_{N,0}(p)(x^N)+xR_{N,1}(p)(x^N)+\cdots +x^{N-1}R_{N,N-1}(p)(x^N) .
$$
For example, when $N=3$
and $p(x)=6x^5+5x^4+4x^3+3x^2+2x+1$ we have $R_{3,0}(p)(x)=4x+1$,
 $R_{3,1}(p)(x)=5x+2$,
 $R_{3,2}(p)(x)=6x+3$. For $N=2$ this corresponds to taking the
 odd ($R_{2,1}$) and even ($R_{2,0}$) parts of the polynomial $p$.

Now, we are ready to establish the following

\proclaim{Theorem}
Suppose $p_n(x)$ ($n=0,1,2, \cdots $) is a sequence of polynomials satisfying
the following $(2N+1)$-term recurrence relation
$$ x^N  p_n(x) = c_{n,0} p_n(x) + \sum_{k=1}^N \left(
  \overline{c_{n,k}} p_{n-k}(x) + c_{n+k,k} p_{n+k}(x) \right) , \tag 2.1$$
where $c_{n,0}$ $(n=0,1,2,\ldots)$ is a real sequence
and $c_{n,k}$ $(n=0,1,2,\ldots)$ are complex sequences
for $k=1,\ldots,N$ with $c_{n,N} \neq 0$ for every $n$
 and with the initial conditions
$p_{k}(x)=0$ for $k <0$ and $p_k$ given polynomials of degree $k$, for $k
=0,\cdots ,N-1$. We define the sequence of matrix polynomials $(P_n)_n$ by
$$
P_n(x)=
  \pmatrix
   R_{N,0}(p_{nN})(x)        &   \cdots & R_{N,N-1}(p_{nN})(x) \\
   R_{N,0}(p_{nN+1})(x)      &   \cdots & R_{N,N-1}(p_{nN+1})(x) \\
      \vdots                 &   \cdots & \vdots   \\
   R_{N,0}(p_{nN+N-1})(x)    &   \cdots & R_{N,N-1}(p_{nN+N-1})(x)
          \endpmatrix, $$
Then this sequence of matrix polynomials is orthonormal on the real line
with respect to a positive definite matrix of measures and satisfies a
matrix three-term recurrence relation. Conversely, suppose
$P_n=(P_{n,m,j})_{m,j=0}^{N-1}$ is a
sequence of orthonormal matrix polynomials or equivalently satisfying a
matrix
three-term recurrence relation (without loss of generality we can assume the
leading coefficient of $P_n$ to be a lower triangular matrix), then the
scalar polynomials defined by
$$
p_{nN+m}(x)=\sum _{j=0}^{N-1}x^jP_{n,m,j}(x^N) ,
\qquad (n\in \N, 0\le m \le N-1),     \tag 2.2
$$
 satisfy a $(2N+1)$-term recurrence
relation of the form \thetag{2.1}.
\endproclaim

\demo{Proof}
The equivalence between $(p_n)_n$ satisfying a $(2N+1)$-term recurrence
relation and $(P_n)_n$ being a sequence of matrix orthonormal polynomials
is a consequence of the definition of the $P_n$ from the $p_k$ (or
conversely, the $p_k$ from the $P_n$)
and the orthogonality condition \thetag{1.4}.

Let us show that the matrix polynomials $(P_n)_n$ satisfy a matrix
three-term recurrence relation, giving explicitly the matrix coefficients which
appear in this recurrence formula. To do that, we consider the $N$-Jacobi
matrix $J$  associated to $(p_n)_n$, which is the $(2N+1)$-banded infinite
Hermitian matrix defined by putting the sequences $(c_{n,l})_n$ which
appear in the recurrence relation on the diagonals of the matrix $J$,
i.e., we define the matrix $J=(j_{n,m})_{n,m \in \N }$ by
$$
j_{n,m}=\cases \overline{c_{n,|n-m|}} & \text{if $0 \leq n-m \leq N$,} \\
   c_{m,|n-m|} &\text{if $0 \leq m-n \leq N$,} \\
   0 &\text{if $|n-m| >N$.} \endcases
$$
Now, we split up this $N$-Jacobi matrix in blocks of dimension $N\times N$,
and then we get the $N\times N$ matrices
$E_n$ and $D_n$ defined by
$$
(E_n)_{i,l}=j_{nN+i, nN+l}=
  \cases \overline{c_{nN+i,|i-l|}} & \text{if $i \geq l$}, \\
         c_{nN+l,|i-l|} & \text{if $i \leq l$,} \endcases
\qquad n \geq 0 , $$
and
$$ \align
(D_n)_{i,l}&=\cases 0 &\text{if $i<l$,} \\ j_{(n-1)N+i,nN+l}  & \text {if
$i\ge l$,} \endcases \\
 &=
\cases 0 &\text {if $i<l$,} \\ c_{nN+l,N+l-i}  &\text {if $i\ge l$,}
\endcases    \qquad n \geq 1,
\endalign
$$
i.e.,
$$  J = \pmatrix  E_0   & D_1   &        &        &        \\
                  D_1^* & E_1   & D_2    &        &        \\
                        & D_2^* & E_2    & D_3    &        \\
                        &       & \ddots & \ddots & \ddots
           \endpmatrix         $$
with
$$   D_n = \pmatrix
    c_{nN,N} &           0 &           0  & \cdots &         0 \\
    c_{nN,N-1}&c_{nN+1,N}   &                0 & \cdots &          0 \\
c_{nN,N-2}    &c_{nN+1,N-1}          & c_{nN+2,N}   & \cdots &          0\\
      \vdots & \ddots       &              & \cdots & \vdots   \\
c_{nN,1}   & c_{nN+1,2}   & c_{nN+2,3}   & \cdots & c_{nN+N-1,N}
          \endpmatrix, $$
and
$$   E_n = \pmatrix
    c_{nN,0}   & c_{nN+1,1}   & c_{nN+2, 2}   & \cdots & c_{nN+N-1, N-1}   \\
 \overline{c_{nN+1,1}} & c_{nN+1,0}  & c_{nN+2,  1} & \cdots & c_{nN+N-1,N-2} \\
\overline{c_{nN+2,2}} & \overline{c_{nN+2,1}}  & c_{nN+2,0}  & \cdots &
c_{nN+N-1, N-3} \\
\vdots    &   \vdots    & \vdots      & \cdots &  \vdots   \\
   \overline{c_{nN+N-1,N-1}} & \overline{c_{nN+N-1,N-2}} &
  \overline{c_{nN+N-1,N-3}} & \cdots & c_{nN+N-1,0} \endpmatrix.  $$
Then  $E_n$ is
an Hermitian matrix, and  the conditions $c_{n,N} \neq 0$ shows that $D_n$ is a
lower triangular matrix
with $\det D_n \neq 0$.

If we compute the following expression
$$
D_{n+1}P_{n+1}(x)+E_nP_n(x)+D_n^*P_{n-1}(x)  \tag 2.3
$$
we find that the entry $(k,m)$ ($0\le k,m\le N-1 $) of this matrix is
equal to
$$
\split \sum _{j=0}^k & c_{(n+1)N+j,N+j-k}R_{N,m}(p_{(n+1)N+j})(x)
\\  & +   \sum_{j=0}^{k-1}\overline{c_{nN+k,|k-j|}}R_{N,m}(p_{nN+j})(x)
 +  \sum_{j=k}^{N-1}c_{nN+j,|k-j|}R_{N,m}(p_{nN+j})(x) \\ &  +
\sum_{j=k}^{N-1} \overline{c_{nN+k,N+k-j}}
     R_{N,m}(p_{(n-1)N+j})(x),  \endsplit  $$
that is,
$$ \multline
c_{nN+k,0}R_{N,m}(p_{nN+k})(x) \\
+  \sum_{l=1}^N \left( \overline{c_{nN+k,l}}R_{N,m}(p_{nN+k-l})(x)+
c_{nN+k+l,l}R_{N,m}(p_{nN+k+l})(x) \right).
\endmultline \tag 2.4
$$
But from the $(2N+1)$-term recurrence relation which  the polynomials
$(p_n)_n$ satisfy,
 it follows that the sequence of polynomials $(R_{N,i}(p_n))_n$
($i=0,\cdots , N-1$) satisfies the following recurrence formula
$$
xR_{N,i}(p_n)(x)=c_{n,0}R_{N,i}(p_n)(x)+\sum_{l=1}^N \left(
\overline{c_{n,l}}R_{N,i}(p_{n-l})(x)+
c_{n+l,l}R_{N,i}(p_{n+l})(x) \right) $$
and so \thetag{2.4} is equal to $xR_{N,m}(p_{nN+k})(x)$. Hence
\thetag{2.3} is equal to $xP_n(x)$ and we have proved that the sequence of
matrix polynomials satisfies the matrix three-term recurrence relation
$$
xP_n(x)=D_{n+1}P_{n+1}(x)+E_nP_n(x)+D_n^*P_{n-1}(x) .
$$
Note that the matrix polynomials obtained in this way always have a leading
coefficient which is a lower triangular matrix and also the matrices
$D_n$ which appear in the recurrence formula are lower triangular.

To prove the converse, suppose $(P_n)_n$ is a sequence of matrix
 polynomials satisfying the recurrence formula
$$
xP_n(x)=D_{n+1}P_{n+1}(x)+E_nP_n(x)+D_n^*P_{n-1}(x) ,
$$
where $\det D_n \neq  0$.
First, we prove that we can assume both the leading coefficients of $P_n$ and
the matrices $D_n$ to be lower triangular matrices, which is necessary
since otherwise the scalar polynomials $p_n$ defined by \thetag{2.2}
could have degree different from $n$.
Orthonormal matrix polynomials
for a given orthogonality measure $\M$ are only determined up to a unitary
factor $U_n$ in the sense that $U_nP_n(x)$ is also orthonormal
with respect to the measure $\M$ whenever $U_nU_n^* = I$. If $Q_n(x)$
is a sequence of orthonormal matrix polynomials satisfying
$$     xQ_n(x) = A_{n+1} Q_{n+1}(x) + B_n Q_n(x) + A_n^* Q_{n-1}(x), $$
then the polynomials $P_n(x) = U_n Q_n(x)$, with $U_nU_n^*=I$, satisfy
$$  xP_n(x) = U_nA_{n+1} U_{n+1}^* P_{n+1}(x) + U_nB_nU_n^* P_n(x)
  + U_nA_n^*U_{n-1}^* P_{n-1}(x) . $$
If the non-singular matrices $A_n$ are not lower triangular, one can
always find  unitary matrices $U_n$ such that
$$    D_n = U_{n-1}A_nU_{n}^* $$
are lower triangular matrices. These unitary matrices are given
recursively
as follows. First we choose a unitary matrix $U_0$ such that $U_0Q_0 = P_0$
is lower triangular. Next we
observe that the matrix
$U_0A_1$ can always be factorized
as $U_0A_1 = D_1U_1$, where $U_1$ is a unitary matrix and $D_1$ is lower
triangular (QR-factorization of Francis and Kublanovskaja). With
this choice of $U_1$ we thus have $D_1=U_0 A_1 U_1^*$. In general, if
$U_0,U_1,\ldots,U_{n-1}$ have been obtained, then we apply the QR-factorization
to $ U_{n-1} A_n$ to find a unitary matrix $U_n$ with
$U_{n-1}A_n=D_nU_n$, with $D_n$ lower triangular.
Therefore any system of orthonormal matrix polynomials can be transformed
to a system of orthonormal matrix polynomials with leading coefficients
which are lower triangular.

Now, since we can assume the leading coefficient of $P_n$ to be  lower
triangular       with non-vanishing determinant,
the sequence of scalar polynomials  $(p_n)_n$ defined  by \thetag{2.2} are such
that dgr($p_n)=n$. And proceeding as before, it is no hard to prove that
they satisfy a $(2N+1)$-term recurrence relation like \thetag{2.1}. The
coefficients
in this recurrence formula can be obtained from the matrices $D_n , E_n$. \qed
\enddemo

For the general case of a $(2N+1)$-term recurrence relation defined by
a polynomial $h$ of degree $N$
$$ h(x) p_n(x) = c_{n,0} p_n(x) + \sum_{k=1}^N \left(
 \overline{c_{n,k}} p_{n-k}(x) + c_{n+k,k} p_{n+k}(x) \right) ,   $$
the theorem works again if
we change the operators $R_{N,m}$ to $R_{h,N,m}$  ($m=0, \cdots , N-1$)
(see \thetag{1.5})
in the definition of the matrix polynomials $P_n$ from the scalar polynomials
$p_n$. Thus, in the definition of the operators $R_{h,N,m}$,
 instead of using the
basis of monomials $\{1,x,x^2,x^3,\ldots\}$ to span the linear space
of polynomials, we will use the basis
$$ \multline
\{ 1,x,\ldots,x^{N-1},h(x),xh(x),\ldots,x^{N-1}h(x),h^2(x),xh^2(x),\ldots\} \\
 = \{ x^n h^m(x): n=0,1,\ldots,N-1,\, m=0,1,2,\ldots \}.
\endmultline \tag 2.5 $$
A polynomial $p$ of degree $nN+m$ $(0 \leq m < N)$ can then be expanded
in this basis as
$$   p(x) =\sum_{i=0}^n \sum_{j=0}^{N-1}  a_{i,j} x^j h^i(x) .$$
Now, the operator $R_{h,N,j}$ takes from $p$ just those terms of the
form $a_{i,j}x^jh^i(x)$ and then, removes the common factor $x^j$ and changes
$h(x)$ to $x$.

Conversely, since now we have
$$
p(x)=R_{h,N,0}(p)(h(x))+xR_{h,N,1}(h(x))+\cdots +x^{N-1}R_{h,N,N-1}(p)(h(x))
$$
we  change \thetag{2.2} in the theorem to
$$
p_{nN+m}(x)=\sum _{j=0}^{N-1}x^jP_{n,m,j}(h(x)),
\qquad (n\in \N, 0\le m \le N-1)
$$
in the definition of the scalar polynomials $p_n$ from the matrix polynomials
$P_n$.

\head 3.  Examples \endhead
\subhead 3.1 Discrete Sobolev orthogonal polynomials \endsubhead

An important class of polynomials satisfying a higher order recurrence
relation is obtained by taking polynomials orthogonal with respect to an
inner product of (discrete) Sobolev type
$$  \< f, g \> = \int f(x) g(x) \, d\mu(x) + \sum_{k=1}^N \lambda_k
        f^{(s_k)}(c_k) q^{(s_k)}(c_k), $$
where $\lambda_k$ are positive real numbers and $c_k$ are real numbers
(which are allowed to coincide). For a discrete Sobolev inner product
$$   \< f,g \> = \int f(x)g(x) \, d\mu(x) + \lambda f^{(r)}(c) g^{(r)}(c) $$
it was shown by Marcell\'an and Ronveaux \cite{16} that the corresponding
monic orthogonal polynomials $q_n(x)$ satisfy a $(2r+3)$-term recurrence
relation
$$  (x-c)^{r+1} q_n(x) = \sum_{j=n-r-1}^{n+r+1} \gamma_{n,j} q_j(x). $$
The orthonormal polynomials then satisfy a recurrence relation of the form
\thetag{1.3} with $h(x)=(x-c)^{r+1}$ and thus these polynomials are also
orthogonal with respect to a positive definite matrix of measures. Evans et al.\
\cite{9}
consider weighted Sobolev inner products
$$   \< f,g \> = \sum_{k=0}^N \int f^{(k)}(x)g^{(k)}(x) \, d\mu_k(x), $$
and show that the existence of a polynomial $h$ for which
$$  \< hp,q \> = \< p,hq \> \tag 3.1 $$
for all polynomials $p,q$ implies that the measures $\mu_k$ $(1 \leq k \leq N)$
are necessarily discrete with support at the zeros of $h$ and that
the orthonormal polynomials satisfy a $(2N+1)$-term recurrence relation
of the form \thetag{1.3}, where $N$ is the degree of the minimal polynomial
$h$ for which \thetag{3.1} holds. It follows that polynomials which are
orthonormal with respect to a discrete  Sobolev inner product correspond
to matrix polynomials orthonormal with respect to a positive definite matrix
of measures.

Consider the discrete Sobolev inner product
$$  \< p,q \> = \int p(x)q(x)\, d\mu(x) + \sum_{i=1}^M \sum_{j=1}^{M_i}
    \lambda_{i,j} p^{(j)}(c_i)q^{(j)}(c_i), \tag 3.2 $$
where $p,q$ are polynomials, $\lambda_{i,j} \geq 0$ and $N = M+\sum_{i=1}^M
M_i$.
Here derivatives are taken at $M$ points $c_i \in {\Bbb R}$ and at the point
$c_i$ the highest derivative is of order $M_i$. Introduce the polynomial
$$   h(x) = \prod_{i=1}^M (x-c_i)^{M_i+1}, $$
then $h$ is of degree $N$ and has its zeros at the points $c_i$ where
the derivatives of the inner product are evaluated. Instead of using the
basis of monomials $\{1,x,x^2,x^3,\ldots\}$ to span the linear space
of polynomials, we will use the basis
$$ \multline
\{ 1,x,\ldots,x^{N-1},h(x),xh(x),\ldots,x^{N-1}h(x),h^2(x),xh^2(x),\ldots\} \\
 = \{ x^n h^m(x): n=0,1,\ldots,N-1,\, m=0,1,2,\ldots \}.
\endmultline \tag 3.3 $$
A polynomial $p$ of degree $Nk+l$ $(0 \leq l < N)$ can then be expanded
in this basis as
$$   p(x) = \sum_{n=0}^{N-1} \sum_{m=0}^k a_{n,m} x^n h^m(x) , $$
where $a_{n,k} = 0$ whenever $n > l$.
Taking the terms in $x^n$ together and putting
$$   R_{h,N,n}(p)(x) = \sum_{m=0}^k a_{n,m} x^m , $$
then for $0 \leq n < N$ each polynomial $R_{h,N,n}(p)$ has degree at
most $k$ (for $n > l$ the degree is less than $k$) and
$$   p(x) = \sum_{n=0}^{N-1} x^n R_{h,N,n}(p)(h(x)) . $$
The polynomial $p$ is thus equivalent (modulo $h$) with the vector polynomial
given by $(R_{h,N,0}(p),R_{h,N,1}(p),\ldots,R_{h,N,N-1}(p))$
and we will write
$$   p \equiv (R_{h,N,0}(p),R_{h,N,1}(p),\ldots,R_{h,N,N-1}(p)) . $$
Observe again that in case $h(x)=x^2$ this
amounts to the decomposition of a polynomial $p$ into its odd and even parts.

For the $j$th derivative of $p$ we then find, using Leibnitz' rule
for the derivative of a product,
$$   p^{(j)}(x) = \sum_{n=j}^{N-1} \sum_{k=0}^j \binom{j}{k}
   n(n-1)\cdots(n-k+1) x^{n-k} \frac{d^{j-k}}{dx^{j-k}} R_{h,N,n}(p)(h(x)) . $$
The derivative of the composite function $R_{h,N,n}(p)(h(x))$ can be evaluated
using the formula of Faa di Bruno, giving
$$ \multline
\frac{d^m}{dx^m} R_{h,N,n}(p)(h(x)) \\
= \sum_{i=0}^m \frac{d^i}{dh^i}
      R_{h,N,n}(p)(h)
  \sum_{a_1,\ldots,a_m} \binom{m}{a_1,a_2,\ldots,a_m} \left( \frac{h'(x)}{1!}
     \right)^{a_1} \cdots \left( \frac{h^{(m)}(x)}{m!}  \right)^{a_m},
\endmultline    $$
where $a_i$ are non-negative integers satisfying $a_1+a_2+\cdots+a_m=i$ and
$a_1+2a_2+\cdots+ma_m=m$. If $m \leq M_i$ and if we evaluate this expression
at $c_i$ then due to the fact that $h'(c_i)=\cdots=h^{(m)}(c_i)=0$ we see that
$$   \left.   \frac{d^m}{dx^m} R_{h,N,n}(p)(h(x))
    \right|_{x=c_i} = 0, \qquad
    1 \leq m \leq M_i. $$
Therefore the only contribution in the expression for $p^{(j)}(c_i)$ when
$1 \leq j \leq M_i$ is when $k=j$, giving
$$  p^{(j)}(c_i) = \sum_{n=0}^{N-1} \frac{n!}{(n-j)!} c_i^{n-j} R_{h,N,n}(p)(0),
  \qquad 1 \leq j \leq M_i. $$
For $p\equiv (R_{h,N,0}(p),\ldots,R_{h,N,N-1}(p))$ and $q\equiv
(R_{h,N,0}(q),\ldots,R_{h,N,N-1}(q))$
the inner product \thetag{3.2} can thus be written as
$$ \align
\< p,q \> &= \int \pmatrix
  R_{h,N,0}(p)(h(x)) &  \ldots & R_{h,N,N-1}(p)(h(x))
\endpmatrix d\M(x) \pmatrix
          R_{h,N,0}(q)(h(x)) \\ R_{h,N,1}(q)(h(x)) \\
 \vdots \\ R_{h,N,N-1}(q)(h(x)) \endpmatrix   \\
&+ \pmatrix R_{h,N,0}(p)(0) &  \ldots &
R_{h,N,N-1}(p)(0) \endpmatrix \L \pmatrix
  R_{h,N,0}(q)(0) \\ R_{h,N,1}(q)(0) \\ \vdots  \\
R_{h,N,N-1}(q)(0) \endpmatrix , \endalign  $$
where $\M$ is the $N\times N$ matrix of measures
$$  d\M(x) = \pmatrix d\mu(x) & x\, d\mu(x) & \cdots & x^{N-1}\,d\mu(x) \\
          x\, d\mu(x) & x^2\, d\mu(x) & \cdots & x^{N}\,d\mu(x) \\
          x^2\, d\mu(x) & x^3\, d\mu(x) & \cdots & x^{N+1}\, d\mu(x) \\
          \vdots        & \vdots        & \cdots & \vdots \\
          x^{N-1}\, d\mu(x) & x^{N}\, d\mu(x) & \cdots & x^{2N-2}\, d\mu(x)
          \endpmatrix  \tag 3.4 $$
and $\L$ is the matrix
$$   \sum_{i=1}^M \sum_{j=1}^{M_i} \lambda_{i,j} \L(i,j)  $$
with $\L(i,j)$ the $N \times N$ matrix
$$  \L(i,j) = \pmatrix  0 \\ \vdots \\ 0 \\
                j!  \\
                \vdots \\ \frac{k!}{(k-j)!} c_i^{k-j} \\
                \vdots \\ \frac{(N-1)!}{(N-1-j)!} c_i^{N-1-j}
                \endpmatrix
            \pmatrix  0  \ldots  0 &
                j!  &
             \ldots & \frac{k!}{(k-j)!} c_i^{k-j} &
                \ldots & \frac{(N-1)!}{(N-1-j)!} c_i^{N-1-j}
                \endpmatrix  $$
with entries
$$   \L_{k,n}(i,j) = \frac{k!\,n!}{(k-j)!\,(n-j)!} c_i^{n+k-2j} ,
\qquad j \leq k,n \leq N-1. $$
It is clear that $\L(i,j)$ is positive definite, hence $\L$ is also
positive definite when all $\lambda_{i,j} \geq 0$.
If $p_n(x)$ $(n=0,1,2,\ldots)$ are the orthonormal polynomials with
respect to the Sobolev inner product \thetag{3.2}, then we can write the
polynomials $p_{kN+l}(x)$ $0 \leq l < N$ using the basis functions
\thetag{3.3}:
$$ p_{kN+l}(x) = \sum_{n=0}^{N-1} x^n R_{h,N,n}(p_{kN+l})(h(x)), \qquad
l=0,1,\ldots,N-1. $$
The matrix polynomials
$$   P_n(x) =
  \pmatrix
   R_{h,N,0}(p_{nN})(x)        &   \cdots & R_{h,N,N-1}(p_{nN})(x) \\
   R_{h,N,0}(p_{nN+1})(x)      &   \cdots & R_{h,N,N-1}(p_{nN+1})(x) \\
      \vdots                 &   \cdots & \vdots   \\
   R_{h,N,0}(p_{nN+N-1})(x)    &   \cdots & R_{h,N,N-1}(p_{nN+N-1})(x)
          \endpmatrix, $$
are then orthonormal with respect to the matrix of measures
$\M(h^{-1})$ to which a mass point at $0$ is added, with weight given
by the matrix $\L$. The matrix of measures $\M(h^{-1})$ is given by
$$ \int F(x) \, d\M(h^{-1}(x)) = \int F(h(x)) \, d\M(x), $$
where $F: {\Bbb R} \to {\Bbb R}^N$ is a vector function such that
$F(h) \in L_1(\M)$.
In this way Sobolev orthogonal polynomials (with a
discrete Sobolev part) can always be expressed as orthogonal matrix
polynomials, where the spectral matrix of measures has a mass point at the
origin. Note that the matrix polynomial $P_k(x)$ is of degree $k$
with a leading coefficient which is a lower triangular matrix because
the degree of $R_{h,N,n}(p_{kN+l})$ is less than $k$ whenever $n > l$.

\subhead 3.2 Perturbation of a measure on the real line by finitely many
 function values \endsubhead

As a second example we consider the inner product
$$  \langle f,g \rangle = \int f(x)g(x) \, d\mu(x)
 + \left( \sum_{k=0}^{N-1} a_kf(c_k) \right) \left( \sum_{k=0}^{N-1} a_kg(c_k)
  \right), \tag 3.5 $$
where $f,g$ are real functions, $\mu$ is a positive measure on the real line
and $a_k,c_k \in {\Bbb R}$ $(k=0,\ldots,N-1)$. When $N=1$ then this is just a
an inner product in $L_2$ with respect to the measure $\mu$ to which
a mass $a_0^2$ is added at the point $c_0$. The case $N=2$ with
$a_0=-a_1$ is basically an inner product involving differences, as introduced
recently by Bavinck \cite{2} \cite{3}.

If we consider the polynomial
$$  h(x) = \prod_{k=0}^{N-1} (x-c_k), $$
then $h$ vanishes at the points $c_k$ where the function values are
evaluated in the inner product, and hence
$$  \langle f,hg \rangle = \int h(x) f(x)g(x) \, d\mu(x) = \langle hf,g \rangle,
$$
which immediately implies that the orthonormal polynomials $(p_n)_n$
for this inner product satisfy a $(2N+1)$-term recurrence formula
of the form \thetag{1.3}. Expanding a polynomial $p$ in the
basis \thetag{3.3} gives
$$  p(x) = \sum_{n=0}^{N-1} x^n R_{h,N,n}(p)(h(x)). $$
Observe that
$$  p(c_k) = \sum_{n=0}^{N-1} c_k^n R_{h,N,n}(p)(0), $$
hence for
 $p\equiv (R_{h,N,0}(p),\ldots,R_{h,N,N-1}(p))$ and $q\equiv
(R_{h,N,0}(q),\ldots,R_{h,N,N-1}(q))$
we have
$$ \multline
\left( \sum_{k=0}^{N-1} a_kf(c_k) \right) \left( \sum_{k=0}^{N-1} a_kg(c_k)
  \right) \\
  = \pmatrix R_{h,N,0}(p)(0) &  \ldots &
R_{h,N,N-1}(p)(0) \endpmatrix \L \pmatrix
  R_{h,N,0}(q)(0) \\ R_{h,N,1}(q)(0) \\ \vdots  \\
R_{h,N,N-1}(q)(0) \endpmatrix ,
\endmultline $$
where
$$  \L = \pmatrix  \sum_{k=0}^{N-1} a_k \\ \sum_{k=0}^{N-1} a_kc_k \\ \vdots \\
                 \sum_{k=0}^{N-1} a_kc_k^{N-1}
                \endpmatrix
            \pmatrix  \sum_{k=0}^{N-1} a_k & \sum_{k=0}^{N-1} a_kc_k & \ldots &
              \sum_{k=0}^{N-1} a_kc_k^{N-1}
                \endpmatrix.  $$
The inner product \thetag{3.5} in the linear space of scalar polynomials
thus corresponds to an inner product in $L_2$ for the measure
$\M(h^{-1})$, with $\M$ given by \thetag{3.4}, to which a Dirac measure at $0$
is added with weight given by the positive definite matrix $\L$.
Moreover, the matrix polynomials
$$   P_n(x) =
  \pmatrix
   R_{h,N,0}(p_{nN})(x)        &   \cdots & R_{h,N,N-1}(p_{nN})(x) \\
   R_{h,N,0}(p_{nN+1})(x)      &   \cdots & R_{h,N,N-1}(p_{nN+1})(x) \\
      \vdots                 &   \cdots & \vdots   \\
   R_{h,N,0}(p_{nN+N-1})(x)    &   \cdots & R_{h,N,N-1}(p_{nN+N-1})(x)
          \endpmatrix $$
are orthonormal with respect to $\M(h^{-1}) + \L \delta_0$.

The special case when one adds to the $L_2(\mu)$ inner product
a part dealing with the $(N-1)$st difference, i.e.,
$$  \sum_{k=0}^{N-1} a_kf(c_k) = \Delta^{N-1} f(c), $$
corresponds to the choice $c_k= c+k\delta$ and
$a_k = (-1)^k \binom{N-1}{k}$. Observe that in this case
$$  \L = \pmatrix 0 \\ 0 \\ \vdots \\ 0 \\ (-\delta)^{N-1}(N-1)! \endpmatrix
        \pmatrix 0 & 0 & \ldots & 0 & (-\delta)^{N-1}(N-1)! \endpmatrix , $$
so that $\L$ contains zeros everywhere, except for the entry in the
lower right corner. This case covers the inner products considered
by Bavinck \cite{2} \cite{3}.
The limiting case where $\delta \to 0$ corresponds to a discrete
Sobolev inner product with a $(N-1)$st derivative at the point $c$.
The general discrete Sobolev inner product \thetag{3.2} can also be
obtained as a limiting case of \thetag{3.5} by letting some of the
$c_k$ coincide and taking appropriate coefficients $a_k$.

\head 4. Krein's theorem revisited \endhead

Suppose that the spectrum $\sigma(J)$ of a (tridiagonal) Jacobi matrix $J$ is
denumerable. If the set of accumulation points of the spectrum $\sigma(J)$
is finite, then a complete characterization of the derived set $\sigma(J)'$ is
given by M.G. Krein \cite{14} \cite{6, Chapter IV, Section 6}:

\proclaim{Krein's theorem}
Suppose $J$ is a bounded Jacobi matrix. Then every accumulation point of
$\sigma(J)$ is a zero of the polynomial $h(x)$ of degree $N$ if and only
if the operator $h(J)$ is compact.
\endproclaim

The relationship between polynomials satisfying a higher order recurrence
relation an orthogonal matrix polynomials enables us to give a short proof of
this result.

\demo{Proof of Krein's theorem}
The infinite matrix $h(J)$ is a banded matrix with band width $2N+1$, and
can thus be written as
$$   h(J) = D + \sum_{j=1}^N (V^*)^j A_j + \sum_{j=1}^N B_j V^j, $$
where $D, A_1, B_1, \ldots, A_N, B_N$ are diagonal matrices and
$V$ is the shift operator that acts as $(V\psi)_n = \psi_{n+1}$.
A diagonal matrix is compact if and only if its entries tend to zero, and
compact operators
on a Hilbert space form a two-sided ideal in the set of bounded operators
on this Hilbert space. Hence the boundedness of the shift operator
implies that $h(J)$ is compact if and only if the entries of $h(J)$ tend to
zero along the diagonals. If we write $h(J)$ as a block Jacobi matrix
$$  h(J) = \pmatrix  E_0   & D_1   &        &        &        \\
                     D_1^* & E_1   & D_2    &        &        \\
                           & D_2^* & E_2    & D_3    &        \\
                           &       & \ddots & \ddots & \ddots
           \endpmatrix         $$
where $E_k, D_k$ are $N\times N$ matrices and $D_k$ are lower triangular,
then the compactness of $h(J)$ implies that $E_k$ and $D_k$ converge
towards the zero matrix. A compact operator has a denumerable spectrum with zero
as the only accumulation point. If $\mu$ is the spectral measure for $J$
(with corresponding orthogonal polynomials $p_n(x)$), then $\M(h^{-1})$
is the matrix of measures for which the matrix polynomials $P_n(x)$
corresponding with $h(J)$ are orthogonal, where
$$  d\M(x) = \pmatrix d\mu(x) & x\, d\mu(x) & \cdots & x^{N-1}\,d\mu(x) \\
          x\, d\mu(x) & x^2\, d\mu(x) & \cdots & x^{N}\,d\mu(x) \\
          x^2\, d\mu(x) & x^3\, d\mu(x) & \cdots & x^{N+1}\, d\mu(x) \\
          \vdots        & \vdots        & \cdots & \vdots \\
          x^{N-1}\, d\mu(x) & x^{N}\, d\mu(x) & \cdots & x^{2N-2}\, d\mu(x)
          \endpmatrix . $$
If $\sigma(h(J))$ is the spectrum of $h(J)$, then it follows that the spectrum
of $J$ is given by $h^{-1}(\sigma(h(J))$. We know that $\sigma(h(J))$ has only
one accumulation point at zero, which therefore corresponds to accumulation
points of $J$ at $h^{-1}(0)$, which are the zeros of the polynomial $h(x)$.

Conversely, if  we know that $\sigma(J)$ is denumerable with accumulation
points at the
zeros of $h$, then the spectrum of $h(J)$ is given by $h(\sigma(J))$ and is
therefore also denumerable with
accumulation points at $h(h^{-1}(0))=0$, which makes
$h(J)$  a compact operator. \qed
\enddemo

\enddocument